\documentclass[12pt]{article}
\usepackage{amsfonts}
\usepackage{amsmath,latexsym,amssymb,amsfonts,amsbsy, amsthm}
\usepackage[usenames]{color}
\usepackage{xspace,colortbl}
\allowdisplaybreaks

\newcommand{\beq}{\begin{equation}}
\newcommand{\eeq}{\end{equation}}
\newcommand{\ben}{\begin{eqnarray}}
\newcommand{\een}{\end{eqnarray}}
\newcommand{\beno}{\begin{eqnarray*}}
\newcommand{\eeno}{\end{eqnarray*}}

\begin{document}

\title{\textbf{Existence of multiple solutions to a class of  nonlinear Schr\"{o}dinger system with external sources terms}
 }
\author{Zexin Qi$^a$,~~ Zhitao Zhang$^{b,}$\thanks{Corresponding author, supported by NSFC 11325107, 11271353, 11331010.} \\
{\small $^a$ College of Mathematics and Information Science,}\\
{\small Henan Normal University, Xinxiang 453007, China}\\
 {\small E-mail: qizedong@126.com}\\
{\small $^b$ Academy of Mathematics and Systems Science,}\\
{\small the Chinese Academy of Sciences, Beijing 100190, China}\\
{\small E-mail: zzt@math.ac.cn} \ }
\date{}
\maketitle
\begin{abstract}
We  study a class of nonlinear schr\"{o}dinger system with  external sources terms  as perturbations in order to obtain existence of multiple solutions,
 this system arises from Bose-Einstein condensates etc..
 As these  external sources   terms are positive functions and small in some sense, we use Nehari manifold to get the existence of a positive ground state solution and a positive  bound state solution.
\end{abstract}


\vskip 0.2cm


\renewcommand{\theequation}{\thesection.\arabic{equation}}
\setcounter{equation}{0}
\section{Introduction}
We are concerned with the following   nonlinear schr\"{o}dinger system with external sources terms
\begin{equation}\label{QZ1.1}
     \left\{
   \begin{array}{ll}
      -\Delta u+\lambda_{1}u=\mu_{1}u^{3}+\beta uv^{2}+f(x), & x\in  \Omega, \\
      -\Delta v+\lambda_{2}v=\mu_{2}v^{3}+\beta u^{2}v+g(x), & x\in \Omega, \\
     u=v=0, & x\in\partial \Omega,
   \end{array}
 \right.
\end{equation}
where $\Omega$ is a bounded domain in~$\mathbb{R}^{N}, N\leq3$,
$\lambda_{1},\lambda_{2},\mu_{1},\mu_{2}$ and $\beta$ are positive
constants, $f(x),g(x)$ are external sources terms. \par (1.1) is a
perturbed version of the following system
\begin{equation}\label{1.2}
     \left\{
   \begin{array}{ll}
      -\Delta u+\lambda_{1}u=\mu_{1}u^{3}+\beta uv^{2}, &  x\in  \Omega, \\
      -\Delta v+\lambda_{2}v=\mu_{2}v^{3}+\beta u^{2}v, &  x\in \Omega, \\
     u=v=0, &  x\in\partial \Omega.
   \end{array}
 \right.
\end{equation}
 \par
 System (1.2) arises from many physical problems, especially in describing some phenomenon in nonlinear optics (\cite{2},\cite{9}). It also models the Hartee-Fock theory for a double
  condensate, i.e., a binary mixture of  Bose-Einstein condensates in two different hyperfine states $|1\rangle$ and $|2\rangle$ (\cite{15}). We refer the reader to
   \cite{3,1,10,ChenZ,CTV2,CTV3,41,42,43,LW,18,NTTV1,NTTV2,34}, and the references therein for  interesting  existence
    of solutions or properties of solutions. The parameters $\mu_i$ and $\beta$ are the intraspecies and interspecies scattering lengths respectively. The sign of the scattering length $\beta$ determines whether the interaction of states $|1\rangle$ and $|2\rangle$ are repulsive or attractive. When $\beta<0$, the interactions of states $|1\rangle$ and $|2\rangle$ are repulsive. In contrast, when $\beta>0$ the interactions of states $|1\rangle$ and $|2\rangle$ are attractive. For atoms of the single state $|j\rangle$, when $\mu_j>0$, the interactions of the single state $|j\rangle$ are attractive.
 \par
  Naturally, people concern nontrivial solutions(solutions with both components nonzero) of the system. In recent years, many interesting works have been devoted to searching ground states and bound states
  for this system, see \cite{3,1,10,ChenZ,LW,18,34} etc. and the references therein.\par
     \textit{\bf A positive ground state solution} we mean a solution of a schr\"{o}dinger system which has the least
      energy among all nonzero solutions, and both of its components are
positive. Note here we call a function positive if it is
  nonnegative and nonzero. A \textit{\bf bound state solution} refers to limited-energy solution.
  As for system (1.2), well-known results indicate the existence of positive ground state is closely related to
   the parameters, see \cite{1},\cite{34} etc. and a remark for the bounded-domain case in \cite{17}.

 In this paper, multiplicity result is established when the perturbations are sufficiently small. If $f(x)$ and $g(x)$ are both positive,
  we  can  find a positive ground state.
  \par
   To be precise, let $S_{4}$ be the best
  Sobolev constant of the embedding: $ H_{0}^{1}(\Omega)\hookrightarrow L^{4}(\Omega)$, then we have\\
 \\
 \textbf{Theorem~1.1}~~\textit{Assume that}~ $f(x),g(x)\in L^{\frac{4}{3}}(\Omega)$,~\textit{both nonzero. Then there exists a positive constant}~$\Lambda =\Lambda(\lambda_{1},\lambda_{2},\mu_{1},\mu_{2}, \beta,S_{4})$, \textit{such that whenever
$\max\{\|f\|_{\frac{4}{3}},\|g\|_{\frac{4}{3}}\}<\Lambda $, system}~(1.1)~\textit{has two nontrivial solutions.~Furthermore,~if~$f$ and $g$~are both positive, system}~(1.1)~
\textit{has one positive ground state solution and one positive bound state solution}.\\

Let $H:=H_{0}^{1}(\Omega)\times H_{0}^{1}(\Omega)$ with  the norm
\begin{center}
$\|(u,v)\|:=\left(\displaystyle\int_{\Omega}(|\nabla u|^{2}+\lambda_{1}u^{2})+\displaystyle\int_{\Omega}(|\nabla v|^{2}+\lambda_{2}v^{2})\right)^{\frac{1}{2}}.$
\end{center}
\par An element ${(}u,v{)}\in H$ is called a weak solution of~(1.1),~if the equality
  $$\displaystyle\int_{\Omega}(
\nabla u \nabla \varphi+\nabla v\nabla \psi-\lambda_{1}u\varphi-\lambda_{2}v\psi-\mu_{1}u^{3}\varphi-\mu_{2}v^{3}\psi$$
$$-\beta uv^{2}\varphi-\beta u^{2}v\psi-f\varphi-g\psi)dx=0~~~~~~~~~~~~$$
holds for all $(\varphi,\psi)\in H$. A weak solution of (1.1) corresponds to a critical point of the following $C^{1}$-functional
\begin{equation}
J(u,v)=\frac{1}{2}\|(u,v)\|^{2}-\frac{1}{4}(\mu_{1}\|u\|_{4}^{4}+\mu_{2}\|v\|_{4}^{4}+2\beta\int_{\Omega}u^{2}v^{2})-\int_\Omega(fu+gv).
\end{equation}
\par Denote the Nehari manifold associated with the functional by
 \begin{center}
 $\mathcal{N}:=\{(u,v)\in H : \langle J^{'}(u,v),(u,v)\rangle=0\}.$
\end{center}
\par It is well-known that all critical points lie in
the~Nehari~manifold, and it is usually  effective to consider
the existence of critical points in this smaller subset of the Sobolev space. For fixed
~$(u,v)\in H \setminus\{(0,0)\}$, denote
\begin{equation}\label{}
\phi(t)=\phi_{(u,v)}(t):=J(tu,tv),~t>0
\end{equation}
the so called fibering map in the direction $(u,v)$.~Such maps are often used to investigate Nehari manifolds for
 various semi-linear problems(\cite{27,Poh,32,40}). Our method is, roughly speaking, to figure out two non-degenerate
  parts of the Nehari manifold, and then consider minimization problems in the two parts
  respectively to obtain two nontrivial solutions, especially to obtain two positive solutions.  Under some other assumptions, there may exist more solutions,
   after we  submitted this manuscript, we thank the referee to tell us
   that existence of infinitely many nontrivial
  solutions of systems related to \eqref{QZ1.1}  are recently obtained  in \cite{YueZ} by different methods under
  different conditions, they   assume  the coefficient matrix is either positive
definite when $N=1,2$ or anti-symmetric when $N=3$. But here for
\eqref{QZ1.1}, the coefficient matrix
  $\left(\begin{array}{lll}
  \mu_1~ \beta\\
  \beta~ \mu_2\end{array}
  \right)$
 usually doesn't satisfy those assumptions, especially the anti-symmetric assumption implies that $\beta=0$ when
 $N=3$, then \eqref{QZ1.1} becomes two independent differential equations without any
 couplings. Our assumption here on $\beta$  can be $\beta>-\sqrt{\mu_1\mu_2}$(notice the remark at the end
 of the paper).
  \par
Our paper is organized as follows. In section 2, we  use fibering
maps to divide the Nehari manifold into three parts. A basic work
related can be found in \cite{Poh},\cite{32}.
 The number $\Lambda$ in Theorem 1.1 is determined to ensure a satisfactory partition. In Section 3 we set up a critical
  point lemma   to consider two minimization problems.
  In the last section, we give the proof of Theorem 1.1.
\section{Partition of the Nehari manifold}
\numberwithin{equation}{section}
 \setcounter{equation}{0}
\par
Denote~$\Phi(u,v):= \langle J^{'}(u,v),(u,v)\rangle$. By (1.3) one has
\begin{center}
$\langle J^{'}(u,v),(u,v)\rangle=\|(u,v)\|^{2}-(\mu_{1}\|u\|_{4}^{4}+\mu_{2}\|v\|_{4}^{4}+2\beta\int_{\Omega}u^{2}v^{2})-\int_\Omega(fu+gv),$
\end{center}
\begin{center}
$\langle \Phi^{'}(u,v),(u,v)\rangle=2\|(u,v)\|^{2}-4(\mu_{1}\|u\|_{4}^{4}+\mu_{2}\|v\|_{4}^{4}+2\beta\int_{\Omega}u^{2}v^{2})-\int_\Omega(fu+gv).$
\end{center}
\par Now we divide the Nehari manifold into three parts:
\begin{center}
   $ \mathcal{N}^{+}:=\{(u,v)\in \mathcal{N} : \langle \Phi^{'}(u,v),(u,v)\rangle>0\};$
\end{center}
\begin{center}
   $ \mathcal{N}^{0}:=\{(u,v)\in \mathcal{N} : \langle \Phi^{'}(u,v),(u,v)\rangle=0\};$
\end{center}
\begin{center}
   $ \mathcal{N}^{-}:=\{(u,v)\in \mathcal{N} : \langle \Phi^{'}(u,v),(u,v)\rangle<0\}.$
\end{center}
\par Obviously, only $\mathcal{N}^{0}$ contains the element $(0,0)$, and it is easy to see $\mathcal{N}^{+}\cup \mathcal{N}^{0}$ and $\mathcal{N}^{-}\cup \mathcal{N}^{0}$ are both closed subsets of $H$.

In order to make an explanation of such partition, for $(u,v)\in H\setminus\{(0,0)\}$,  let us consider the fibering map defined in (1.4).
Since
\begin{center}
 $\phi'(t)=\dfrac{1}{t}\langle J'(tu,tv),(tu,tv)\rangle,$
\end{center}
we know $(u,v)\in \mathcal{N}$ if and only if $\phi'(1)=0$. Moreover, $(tu,tv)\in \mathcal{N}$ with $t>0$ if and only if $\phi'(t)=0$. \par
Thus for a fixed direction $(u,v)\in H \setminus\{(0,0)\}$, we can obtain all elements on the Nehari manifold which lie in this direction if we can find all stationary points of the fibering map. As a result, one would obtain all nonzero elements of the Nehari manifold if one could find stationary points of fibering maps in all directions. We remark that the number of roots of the equation $\phi'(t)=0$ does not depend on the norm of $(u,v)$, once this direction is fixed. Indeed, for $\delta>0$ we have
\begin{center}
$\phi_{(\delta u,\delta v)}^{'}\left(\dfrac{t}{\delta}\right)=\delta\phi_{(u,v)}^{'}(t)$.
\end{center}
 \par Thus
\begin{center}
 $\phi_{(u,v)}^{'}(t)=0,~t>0\Longleftrightarrow\phi_{(\delta u,\delta v)}^{'}\left(\dfrac{t}{\delta}\right)=0,~t>0.$
\end{center}
\par Furthermore,
\begin{center}
$\phi_{(\delta u,\delta v)}^{''}\left(\dfrac{t}{\delta}\right)=\delta^{2}\phi_{(u,v)}^{''}(t)$,
\end{center}
which means it is also irrelevant with the norm of $(u,v)$ if we consider the sign of second derivative of the fibering map at its stationary points.
Moreover, stationary points can be classified into three types, namely local minimum, local maximum and turning point, according to the sign of second derivative of the fibering map at these points.\par
 We now disclose the relationship between such classification and the former partition of the Nehari manifold. By direct calculation, we get
\begin{center}
 $\phi^{'}(t)=\dfrac{1}{t} \langle J^{'}(tu,tv),(tu,tv)\rangle=\dfrac{1}{t}\Phi(tu,tv),$
\end{center}
\begin{center}
 $\phi^{''}(t)=\dfrac{1}{t^{2}}[\langle\Phi^{'}(tu,tv),(tu,tv)\rangle-\Phi(tu,tv)].$
\end{center}
\par Thus if~$\phi^{'}(t)=0$,~then $(tu,tv)\in \mathcal{N}$, and $\Phi(tu,tv)=0$, which yields
\begin{center}
 $\phi^{''}(t)=\dfrac{1}{t^{2}} \langle\Phi^{'}(tu,tv),(tu,tv)\rangle.$
\end{center}
\par Now it is easy to check:
\begin{center}
$t(u,v)\in \mathcal{N}^{+},~t>0 \Longleftrightarrow \phi^{'}(t)=0, \phi^{''}(t)>0;$
\end{center}
\begin{center}
$t(u,v)\in \mathcal{N}^{0},~t>0 \Longleftrightarrow \phi^{'}(t)=0, \phi^{''}(t)=0;$
\end{center}
\begin{center}
$t(u,v)\in \mathcal{N}^{-},~t>0 \Longleftrightarrow \phi^{'}(t)=0, \phi^{''}(t)<0.$
\end{center}
\par We can explore the Nehari manifold through fibering maps. In
fact,  we give the following important lemma to show
when  the degenerate part of the Nehari manifold is clear and
simple. \par
To simplify the calculation, let us introduce some notations
that will be used repeatedly in the rest. For $(u,v)\in H$, define
\begin{equation}\label{}
A=A(u,v):=\mu_{1}\|u\|_{4}^{4}+\mu_{2}\|v\|_{4}^{4}+2\beta\int_{\Omega}u^{2}v^{2},
\end{equation}
\begin{equation}\label{}
B=B(u,v):=\int_\Omega(fu+gv).
\end{equation}
\textbf{Lemma~2.1}~~\textit{Suppose that $f(x),g(x)\in L^{\frac{4}{3}}(\Omega)$,~both nonzero,~then there exists a positive constant~$\Lambda =\Lambda(\lambda_{1},\lambda_{2},\mu_{1},\mu_{2}, \beta,S_{4})$
 such that  $\mathcal{N}^{0}=\{(0,0)\}$ when
$\max\{\|f\|_{\frac{4}{3}},\|g\|_{\frac{4}{3}}\}<\Lambda$.\\}
\\
\textbf{Proof.}~~By the analysis above, we only need to prove for~each $(u,v)\in H $ with $\|(u,v)\|=1$,~the fibering map~$\phi(t)=\phi_{(u,v)}(t)$~has no stationary point that is a turning point. By the notations  (2.1) and (2.2), we can write
\begin{align}\label{}
    \phi(t)&=J(tu,tv)\notag\\
           &=\dfrac{t^{2}}{2}\|(u,v)\|^{2}-\dfrac{t^{4}}{4}(\mu_{1}\|u\|_{4}^{4}+\mu_{2}\|v\|_{4}^{4}+2\beta\displaystyle\int_{\Omega}u^{2}v^{2})-t\int_\Omega(fu+gv)\notag\\
           &=\dfrac{t^{2}}{2}-\dfrac{At^{4}}{4}-Bt.\notag
\end{align}
\par Thus
\begin{center}
   $\phi'(t)=t-At^{3}-B$,\;\;\;\;\;\;\;\;$\phi''(t)=1-3At^{2}.$
\end{center}
\par Notice~$A=A(u,v)=\mu_{1}\|u\|_{4}^{4}+\mu_{2}\|v\|_{4}^{4}+2\beta\int_{\Omega}u^{2}v^{2}>0$.~Define $\psi(t)=t-At^{3},~t>0$, then $\phi'(t)=0\Longleftrightarrow \psi(t)=B$.
\par
Let us consider the graph of $\psi(t)$: $\psi''(t)=-6At<0$,~so~$\psi(t)$~is strictly concave;~$\psi'(t)=0\Longleftrightarrow t=\sqrt{\frac{1}{3A}}$,~thus~$\psi(t)$~takes its maximum at~$t=\sqrt{\frac{1}{3A}}$, the value of which is~$\frac{2}{3}\sqrt{\frac{1}{3A}}$. Also, $\lim_{t\rightarrow0^{+}} \psi(t)=0$,~$ \psi(\infty)=-\infty$. All the above implies if~$0<B<\frac{2}{3}\sqrt{\frac{1}{3A}}$, the equation $\phi'(t)=0$~has exactly two roots $t_{1}$,~$t_{2}$ satisfying $0<t_{1}<\sqrt{\frac{1}{3A}}<t_{2}$; if~$B\leq0$ the equation $\phi'(t)=0$ has exactly one root denoted $t_{2}^{'}$, which satisfies $\sqrt{\frac{1}{3A}}<t_{2}^{'}$. Since~$\phi''(t)=1-3At^{2}$,~considering the above two cases we have $\phi''(t_{1})>0$,~$\phi''(t_{2})<0$, and $\phi''(t_{2}^{'})<0$. So when $0<B<\frac{2}{3}\sqrt{\frac{1}{3A}}$, one has $t_{1}(u,v)\in \mathcal{N}^{+}$,~$t_{2}(u,v)\in \mathcal{N}^{-}$;~
when $B\leq0$ one has $t_{2}^{'}(u,v)\in \mathcal{N}^{-}$. Since $f(x),g(x)$ are both nonzero, it is easy to check that the sets $\{(u,v)\in H :\|(u,v)\|=1,0<B<\frac{2}{3}\sqrt{\frac{1}{3A}}\}\not=\emptyset$ and $\{(u,v)\in H :\|(u,v)\|=1,B\leq0 \}\not=\emptyset$, so $\mathcal{N}^{+}\not=\emptyset$ and $\mathcal{N}^{-}\not=\emptyset$. To finish the proof, we are in a position to determine a number $\Lambda$ such that whenever
$\max\{\|f\|_{\frac{4}{3}},\|g\|_{\frac{4}{3}}\}<\Lambda$ we have
\begin{equation}\label{}
     B<\frac{2}{3}\sqrt{\frac{1}{3A}}.
\end{equation}
\par In fact, since $(u,v)$ lies on the unit sphere of $H$, by Sobolev inequality and~H\"{o}lder~inequality one obtains an upper bound for
~$A$, which yields the existence of a positive constant $\alpha=\alpha(\lambda_{1},\lambda_{2},\mu_{1},\mu_{2},\beta,S_{4})$ such that
\begin{equation}\label{}
        0<\alpha\leq\frac{2}{3}\sqrt{\frac{1}{3\sup\limits_{\|(u,v)\|=1}A(u,v)}}.
\end{equation}
\par On the other hand, by Sobolev inequality and~H\"{o}lder inequality, we obtain that
 \begin{center}
  $B=\displaystyle\int_\Omega(fu+gv)
   \leq \|f\|_{\frac{4}{3}}\|u\|_{4}+\|g\|_{\frac{4}{3}}\|v\|_{4}
    \leq \sqrt{2}S_{4}\max\{\|f\|_{\frac{4}{3}},\|g\|_{\frac{4}{3}}\}\|(u,v)\|.$
\end{center}
\par That is,
\begin{equation}\label{}
B\leq \sqrt{2}S_{4}\max\{\|f\|_{\frac{4}{3}},\|g\|_{\frac{4}{3}}\}.
\end{equation}
\par Now take $\Lambda:=\Lambda(\lambda_{1},\lambda_{2},\mu_{1},\mu_{2},\beta,S_{4}) =\frac{\alpha}{\sqrt{2}S_{4}}$,
by (2.4),(2.5) we know when
$\max\{\|f\|_{\frac{4}{3}},\|g\|_{\frac{4}{3}}\}<\Lambda$,~(2.3) holds.\hfill{$\square$}

\section{Reduction to minimization problems}
\numberwithin{equation}{section}
 \setcounter{equation}{0}
Here we reduce our discussion into two minimization problems through
a critical point lemma. As to constraint minimization problems, we
refer the reader to \cite{8}. Let $X$,$Y$ be real Banach spaces,
$U\subset X$ be an open set. Suppose that $f:U\rightarrow
\mathbb{R}^{1}$, $g:U\rightarrow Y$ are $C^{1}$ mappings in this
section. Let
\begin{center}
   $ M=\{x\in U:g(x)=\theta\}.$
\end{center}
\par For the following minimization problem
\begin{equation}\label{}
\min_{x\in M}f(x),
  \end{equation}
we have\\
\\
\textbf{Lemma~3.1}~(\cite{8},~Theorem 4.1.1)~~\textit{Suppose that~$x_{0}\in M$~solves}~(3.1),~\textit{and that~$Im~g'(x_{0})$~is closed.~Then there exists~$(\lambda,y^{\ast})\in \mathbb{R}^{1}\times Y^{\ast}$~such that~$(\lambda,y^{\ast})\neq (0,\theta)$,~and
\begin{center}
$\lambda f'(x_{0})+ g'(x_{0})^{\ast}y^{\ast}=0.$
\end{center}
Furthermore,~if~$Im~g'(x_{0})=Y$,~then~$\lambda\neq 0$.}\\

Let us introduce two minimization problems
\begin{equation}\label{}
 \theta^{+}= \inf_{(u,v)\in \mathcal{N}^{+}}J(u,v);
  \end{equation}
\begin{equation}\label{}
 \theta^{-}= \inf_{(u,v)\in \mathcal{N}^{-}}J(u,v).
  \end{equation}
\textbf{Lemma ~3.2}~~\textit{If~$(u_{1},v_{1})$~solves}~(3.2),~\textit{then~$(u_{1},v_{1})$~is a nontrivial weak solution of}~(1.1);~
\textit{if~$(u_{2},v_{2})$~solves}~(3.3),~\textit{then~$(u_{2},v_{2})$~is a nontrivial weak solution of}~(1.1).\\
\\
\textbf{Proof.}~~We prove the first assertion, since the second is similar. \par
Recall~$\Phi(u,v)= \langle J^{'}(u,v),(u,v)\rangle$~and
\begin{center}
   $ \mathcal{N}^{+}=\{(u,v)\in \mathcal{N} : \langle \Phi^{'}(u,v),(u,v)\rangle>0\}.$
\end{center}
\par Let~$U=\{(u,v)\in H :\langle \Phi^{'}(u,v),(u,v)\rangle>0\}$,~and rewrite
\begin{center}
   $ \mathcal{N}^{+}=\{(u,v)\in U:\Phi(u,v)=0 \}.$
\end{center}
\par Now we use Lemma 3.1 to consider problem (3.2) by taking
$X=H,~Y=\mathbb{R}^{1},~U=\{(u,v)\in H :\langle
\Phi^{'}(u,v),(u,v)\rangle>0\},~f=J,~g=\Phi,~M=\mathcal{N}^{+}$.
Then there exists a real pair~$(\lambda_{1},\lambda_{2})\neq
(0,0)$,~satisfying
\begin{center}
    $\lambda_{1}J'(u_{1},v_{1})=\lambda_{2}\Phi'(u_{1},v_{1}).$
\end{center}
\par Since $\langle \Phi^{'}(u_{1},v_{1}),(u_{1},v_{1})\rangle>0$, we
have~$\Phi^{'}(u_{1},v_{1})\neq 0$,~again by Lemma 3.1,~one may
assume~$\lambda_{1}\neq0$.~In other words, there exists a real
number~$\mu$ such that
\begin{equation}\label{}
    J'(u_{1},v_{1})=\mu\Phi'(u_{1},v_{1}).
\end{equation}
\par Since
\begin{center}
   $\mu \langle \Phi^{'}(u_{1},v_{1}),(u_{1},v_{1})\rangle=\langle J^{'}(u_{1},v_{1}),(u_{1},v_{1})\rangle=0$
\end{center}
and~$\langle \Phi^{'}(u_{1},v_{1}),(u_{1},v_{1})\rangle>0$,~we
obtain~$\mu =0$. By (3.4)~we get~$J'(u_{1},v_{1})=0$.~
Obviously~$(u_{1},v_{1})\neq (0,0)$, noticing that (1.1)~doesn't
admit semi-trivial (nonzero but one component being trivial)
solutions, we know that ~$(u_{1},v_{1})$~is a nontrivial weak
solution of (1.1).\hfill{$\square$}
\section{Proof of the main result}
\numberwithin{equation}{section}
 \setcounter{equation}{0}
When $\mathcal{N}^{0}=\{(0,0)\}$, we show the solvability of the two
problems (3.2) and (3.3), and give the proof of Theorem 1.1.
The following lemma implies that~$\theta^{+}$~and~$\theta^{-}$~are both finite,~and any minimizing sequence for (3.2) or (3.3) are bounded.\\
\\
\textbf{Lemma~4.1}~~\textit{Assume that $f(x),g(x)\in L^{\frac{4}{3}}(\Omega)$,~then~ $J$~ is coercive and bounded from below on~$\mathcal{N}$~}(\textit{thus on~$\mathcal{N}^{+}$~and~$\mathcal{N}^{-}$}).\\
\\
\textbf{Proof.}~~Let~$(u,v)\in \mathcal{N}$, from the definition of the~Nehari~manifold we have
\begin{center}
   $
   \|(u,v)\|^{2}-\displaystyle\int_\Omega(fu+gv)=(\mu_{1}\|u\|_{4}^{4}+\mu_{2}\|v\|_{4}^{4}+2\beta\displaystyle\int_{\Omega}u^{2}v^{2})$,
\end{center}
this equality together with~(1.3)~yield
\begin{equation}
J(u,v)=\frac{1}{4}\|(u,v)\|^{2}-\frac{3}{4}\int_\Omega(fu+gv).
\end{equation}
\par Note that
\begin{equation}\label{}
   \displaystyle\int_\Omega(fu+gv)
   \leq \|f\|_{\frac{4}{3}}\|u\|_{4}+\|g\|_{\frac{4}{3}}\|v\|_{4}
    \leq \sqrt{2}S_{4}\max\{\|f\|_{\frac{4}{3}},\|g\|_{\frac{4}{3}}\}\|(u,v)\|,
  \end{equation}
combining (4.1) and (4.2) we obtain
\begin{equation}
J(u,v)\geq \frac{1}{4}\|(u,v)\|^{2}-\frac{3\sqrt{2}}{4}S_{4}\max\{\|f\|_{\frac{4}{3}},\|g\|_{\frac{4}{3}}\}\|(u,v)\|.
\end{equation}
\par Since the right hand side of the inequality (4.3) is a quadratic function of $\|(u,v)\|$, it is easy to know that $J$ is coercive and bounded from below on~$\mathcal{N}$.\hfill{$\square$}\\

From now on, $\Lambda$ refers to the number in Lemma 2.1. Although
the lemma below is valid under weaker conditions, we would assume
$f(x),g(x)\in L^{\frac{4}{3}}(\Omega)$, both nonzero, and
$\max\{\|f\|_{\frac{4}{3}},\|g\|_{\frac{4}{3}}\}<\Lambda$ in order
to ensure $\mathcal{N}^{0}=\{(0,0)\}$ by Lemma 2.1, because our discussion bases on a good partition of the
Nehari manifold. We  will often use the proof
of Lemma 2.1 to analyze some properties of the manifold as well as
associated fibering maps.

 For the first minimization problem, we establish:\\
\\
\textbf{Lemma~4.2}~~\textit{Assume that $f(x),g(x)\in L^{\frac{4}{3}}(\Omega)$, both nonzero, and that $\max\{\|f\|_{\frac{4}{3}},\|g\|_{\frac{4}{3}}\}<\Lambda$. Then $\theta^{+} < 0$}.\\
\\
\textbf{Proof.}~~By the proof of Lemma 2.1, we may take one $(u,v)\in H $,~$\|(u,v)\|=1$ such that $0<B<\frac{2}{3}\sqrt{\frac{1}{3A}}$. Now the derivative  $\phi'(t)$ of the
 fibering map in this direction has exactly two positive zeros: $t_{1}$,~$t_{2}$, satisfying~$0<t_{1}<\sqrt{\frac{1}{3A}}<t_{2}$, and $t_{1}\cdot(u,v)\in \mathcal{N}^{+}$.
By $\phi'(t)=t-At^{3}-B$ we have $\lim_{t\rightarrow0^{+}}
\phi'(t)=-B<0$, and $\phi''(t)>0$, $\forall
t\in(0,\sqrt{\frac{1}{3A}})$. Since $\phi'(t_{1})=0$, we know that
~$\phi(t_{1})<\lim_{t\rightarrow0^{+}} \phi(t)=0$. Since
~$\phi(t_{1})=J(t_{1}u,t_{1}v)\geq\theta^{+}$, we obtain that $\theta^{+} < 0$.\hfill{$\square$}\\

For the second minimization problem we need the following lemma, from which we know $\mathcal{N}^{-}$~stays away from the origin.\\
\\
\textbf{Lemma~4.3}~~\textit{Assume that $f(x),g(x)\in L^{\frac{4}{3}}(\Omega)$, both nonzero, and that $\max\{\|f\|_{\frac{4}{3}},\|g\|_{\frac{4}{3}}\}<\Lambda$. Then $\mathcal{N}^{-}$ is closed}.\\
\\
\textbf{Proof.} Under the  assumptions one has
\begin{center}
    $cl~({{\mathcal{N}^{-}}})\subset\mathcal{N}^{-}\cup\mathcal{N}^{0}=\mathcal{N}^{-}\cup \{(0,0)\},$
\end{center}
where we denote $cl~({{\mathcal{N}^{-}}})$ the closure of
$\mathcal{N}^{-}$. Thus we only need to prove $(0,0) \not\in
cl~({{\mathcal{N}^{-}}})$, which is equivalent to prove that
\begin{center}
    ~$dist((0,0),\mathcal{N}^{-})>0.$
\end{center}
\par Take~$(u,v)\in\mathcal{N}^{-}$,~and denote
\begin{center}
$(\widetilde{u},\widetilde{v})=\dfrac{(u,v)}{\|(u,v)\|},$
\end{center}
then~$\|(\widetilde{u},\widetilde{v})\|=1$. Under the assumptions,
by the proof of Lemma 2.1, we obtain that
$A:=A(\widetilde{u},\widetilde{v}),~B:=B(\widetilde{u},\widetilde{v})$
satisfy $B<\frac{2}{3}\sqrt{\frac{1}{3A}}$.~Furthermore,
if~$0<B<\frac{2}{3}\sqrt{\frac{1}{3A}}$,~the
equation~$\phi'_{(\widetilde{u},\widetilde{v})}(t)=0$~has exactly
two roots also denoted by~$t_{1}$, $t_{2}$,~which
satisfy~$t_{2}\cdot(\widetilde{u},\widetilde{v})\in
\mathcal{N}^{-}$,~$t_{1}\cdot(\widetilde{u},\widetilde{v})\in
\mathcal{N}^{+}$, we have $t_{2}\cdot(\widetilde{u},\widetilde{v})=
(u,v)$,~so~$t_{2}= \|(u,v)\|$.~If~$B\leq 0$,~then the equation
~$\phi'_{(\widetilde{u},\widetilde{v})}(t)=0$~has exactly one
root~still denoted by $\widetilde{t}_{2}$,~thus we get
~$\widetilde{t}_{2}\cdot(\widetilde{u},\widetilde{v})= (u,v)$, then
$\widetilde{t}_{2}= \|(u,v)\|$. Since~${t_{2}}>\sqrt{\frac{1}{3A}}$,
$\widetilde{t}_{2}>\sqrt{\frac{1}{3A}}$ in the proof of Lemma
2.1,~so no matter which of the above two cases happens,~ we always
obtain~$\|(u,v)\|>\sqrt{\frac{1}{3A}}$. Noticing that $A$ is bounded
from above,~ we know that there exists~$\tau>0$~such that
\begin{center}
    $\|(u,v)\|>\tau.$
\end{center}
\par We obtain
\begin{center}
    $dist((0,0),\mathcal{N}^{-})=\inf\limits_{(u,v)\in\mathcal{N}^{-}}\|(u,v)\|\geq\tau>0,$
\end{center}
which completes the proof.\hfill{$\square$}\\

In order to abstract a $(PS)_{\theta^{+}}$~sequence from the minimizing sequence for problem (3.2), we use the idea of \cite{32} to obtain the following lemma. \\
\\
\textbf{Lemma~4.4}~~\textit{Assume that $f(x),g(x)\in
L^{\frac{4}{3}}(\Omega)$, both nonzero. Then for $(u,v)\in
\mathcal{N}^{+}$, there exists~$\epsilon=\epsilon(u,v)>0$~and a
differentiable
function~$\xi^{+}:B_{\epsilon}(0,0)\rightarrow\mathbb{R}_{+}:=(0,+\infty)$
such that}
\begin{itemize}
  \item $\xi^{+}(0,0)=1;$
  \item $\xi^{+}(w,z)(u-w,v-z)\in \mathcal{N}^{+},\forall(w,z)\in B_{\epsilon}(0,0);$
  \item $ \langle (\xi^{+})'(0,0),(w,z)\rangle=[\|(u,v)\|^{2}-3A(u,v)]^{-1}[2\displaystyle\int_{\Omega}(\nabla u\nabla w+\lambda_{1}uw+\nabla v\nabla z+\lambda_{2}vz)-4\displaystyle\int_{\Omega}(\mu_{1}u^{3}w+\mu_{2}v^{3}z+\beta uv^{2}w+\beta u^{2}vz)-\displaystyle\int_{\Omega}(fw+gz)].$
\end{itemize}
\textbf{Proof.}~~Define a $C^{1}$-mapping~$F:\mathbb{R}_{+}\times H \rightarrow \mathbb{R}$ as follows:
\begin{center}
    $F(t,(w,z))=t\|(u-w,v-z)\|^{2}-t^{3}A(u-w,v-z)-\displaystyle\int_\Omega[f\cdot(u-w)+g\cdot(v-z)].$
\end{center}
\par
 By (2.1), we know that
$A(u-w,v-z)=\mu_{1}\|u-w\|_{4}^{4}+\mu_{2}\|v-z\|_{4}^{4}+2\beta\displaystyle\int_{\Omega}(u-w)^{2}(v-z)^{2}$.
Since $(u,v)\in \mathcal{N}^{+}$~we have
\begin{center}
    $F(1,(0,0))=0.$
\end{center}
\par Consider the fibering map~$\phi(t)=\phi_{(u,v)}(t):=J(tu,tv)$,~since
\begin{center}
   $F(t,(0,0))=t\|(u,v)\|^{2}-t^{3}(\mu_{1}\|u\|_{4}^{4}+\mu_{2}\|v\|_{4}^{4}+2\beta\displaystyle\int_{\Omega}u^{2}v^{2})-\displaystyle\int_\Omega(fu+gv),$
\end{center}
we have
\begin{center}
   $F(t,(0,0))=\phi'(t).$
\end{center}
\par Since $(u,v)\in \mathcal{N}^{+}$, we get $\phi''(1)>0$,~thus
\begin{center}
   $\dfrac{\partial F}{\partial t}(1,(0,0))=\phi''(1)>0.$
\end{center}
\par We apply the implicit function theorem at point (1,(0,0)) to obtain the existence of~$\epsilon=\epsilon(u,v)>0$~and differentiable
function~$\xi^{+}(u,v)(i.e.,~t(u,v)):B_{\epsilon}(0,0)\rightarrow\mathbb{R}_{+}$
such that
\begin{itemize}
  \item $\xi^{+}(0,0)=1;$
  \item $\xi^{+}(w,z)\cdot(u-w,v-z)\in \mathcal{N},~\forall(w,z)\in B_{\epsilon}(0,0).$
\end{itemize}
\par Besides, we obtain the third conclusion of this lemma by calculation.~To finish the second one, we only need to choose $\epsilon=\epsilon(u,v)>0$ small,~such that~$\xi^{+}(w,z)(u-w,v-z)\in \mathcal{N^{+}},~\forall(w,z)\in B_{\epsilon}(0,0)$.~Indeed,~since
$\mathcal{N}^{-}\cup \mathcal{N}^{0}$~is
closed,~$dist((u,v),\mathcal{N}^{-}\cup \mathcal{N}^{0})>0$. Since
~$\xi^{+}(w,z)\cdot(u-w,v-z)$~is continuous with~respect to
$(w,z)$,~when~$\epsilon=\epsilon(u,v)>0$~is small enough, we know
\begin{center}
   $ \|\xi^{+}(w,z)\cdot(u-w,v-z)-(u,v)\| < \dfrac{1}{2}dist((u,v),\mathcal{N}^{-}\cup \mathcal{N}^{0}),\forall(w,z)\in B_{\epsilon}(0,0).$
\end{center}
That is,~$\xi^{+}(w,z)\cdot(u-w,v-z)$~does not belong to~$
\mathcal{N}^{-}\cup \mathcal{N}^{0}$.~Thus
~$\xi^{+}(w,z)\cdot(u-w,v-z)\in \mathcal{N^{+}}$ and our proof is completed.\hfill{$\square$}\\

Similarly, we can establish the following lemma, which will be used to abstract a $(PS)_{\theta^{-}}$~sequence from the minimizing sequence for problem (3.3).\\
\\
\textbf{Lemma~4.5}~~\textit{Assume that $f(x),g(x)\in
L^{\frac{4}{3}}(\Omega)$, both nonzero. Then for $(u,v)\in
\mathcal{N}^{-}$, there exists~$\epsilon=\epsilon(u,v)>0$~and a
differentiable
function~$\xi^{-}:B_{\epsilon}(0,0)\rightarrow\mathbb{R}_{+}$ such
that}
\begin{itemize}
  \item $\xi^{-}(0,0)=1;$
  \item $\xi^{-}(w,z)(u-w,v-z)\in \mathcal{N}^{-},\forall(w,z)\in B_{\epsilon}(0,0);$
  \item $ \langle (\xi^{-})'(0,0),(w,z)\rangle=[\|(u,v)\|^{2}-3A(u,v)]^{-1}[2\displaystyle\int_{\Omega}(\nabla u\nabla w+\lambda_{1}uw+\nabla v\nabla z+\lambda_{2}vz)-4\displaystyle\int_{\Omega}(\mu_{1}u^{3}w+\mu_{2}v^{3}z+\beta uv^{2}w+\beta u^{2}vz)-\displaystyle\int_{\Omega}(fw+gz)].$
\end{itemize}

We are in a position to give:\\
\\
\textbf{Lemma~4.6}~~\textit{Assume that $f(x),g(x)\in
L^{\frac{4}{3}}(\Omega)$, both nonzero, and that
$\max\{\|f\|_{\frac{4}{3}},\|g\|_{\frac{4}{3}}\}<\Lambda$. Then
there exists a sequence~$\{(u_{n},v_{n})\}\subset
\mathcal{N}^{+}$~such that ($n \rightarrow \infty$)}:
\begin{enumerate}
  \item $J(u_{n},v_{n})\rightarrow \theta^{+};$
  \item $J'(u_{n},v_{n})\rightarrow 0.$
\end{enumerate}
\textbf{Proof.}~~Notice~$\mathcal{N}^{+}\cup \{(0,0)\}$~is closed
in~$H$,~we use the Ekeland's variational principle(\cite{14}) on
$\mathcal{N}^{+}\cup \{(0,0)\}$ to obtain a minimizing sequence
$\{(u_{n},v_{n})\}\subset \mathcal{N}^{+}\cup \{(0,0)\}$ such that
\begin{description}
  \item[(a)] $J(u_{n},v_{n})< \inf\limits_{(u,v)\in \mathcal{N}^{+}\bigcup \{(0,0)\}}J(u,v)+\dfrac{1}{n}$;
  \item[(b)] $J(w,z)\geq J(u_{n},v_{n}) -\dfrac{1}{n}\|(w-u_{n},z-v_{n})\|,\;\;\forall (w,z)\in \mathcal{N}^{+}\cup \{(0,0)\}$.
\end{description}
\par
By Lemma~4.2, we know $\theta^{+} < 0$; since $ J(0,0)=0$, we get
that $$ \inf\limits_{(u,v)\in \mathcal{N}^{+}\bigcup
\{(0,0)\}}J(u,v)=\theta^{+}.$$ \par
Thus $J(u_{n},v_{n})\rightarrow
\theta^{+}$,   we may assume $\{(u_{n},v_{n})\}\subset
\mathcal{N}^{+}$, then the first assertion holds.

For the second assertion, Firstly we have   that
$\inf\limits_n\|(u_{n},v_{n})\|\geq m>0$, where $m$ is a constant. Indeed, if
not, then $J(u_{n},v_{n})$~would converge to~zero. Moreover, by
Lemma~4.1 we know that $J$~is coercive~on $\mathcal{N}^{+}$,~ then
$\{\|(u_{n},v_{n})\|\}$ is bounded.  i.e., $\exists M>0$ such that
\begin{equation}\label{}
    0<m\leq \|(u_{n},v_{n})\| \leq M.
\end{equation}
\par Now by contradiction we assume $\|J'(u_{n},v_{n})\|\geq C>0$ as $n$
is large, otherwise we may extract a subsequence to get the
conclusion.

Now let us take~$(u,v)=(u_{n},v_{n})$ in Lemma 4.4,~and define a differentiable function~$\xi_{n}^{+}:(-\epsilon,\epsilon)\rightarrow \mathbb{R}_{+}$~as:
\begin{center}
   $\xi_{n}^{+}(\delta ):= \xi^{+}\left(\dfrac{\delta J'(u_{n},v_{n})}{\|J'(u_{n},v_{n})\|}\right),$
\end{center}
\par then by  Lemma 4.4, we know that $\xi_{n}^{+}(0)=
\xi^{+}(0,0)=1$,~and for~$ \delta \in (-\epsilon,\epsilon)$~we have
\begin{center}
   $(w,z)_{\delta}:=\xi_{n}^{+}(\delta )\cdot[(u_{n},v_{n})-\dfrac{\delta J'(u_{n},v_{n})}{\|J'(u_{n},v_{n})\|}]\in\mathcal{N}^{+}.$
\end{center}
\par
Since~$(u_{n},v_{n})$~satisfies (b),~one has
\begin{equation}\label{}
J(u_{n},v_{n})-J((w,z)_{\delta})\leq \dfrac{1}{n}\|(w,z)_{\delta}-(u_{n},v_{n})\|.
\end{equation}
\par Expanding the left hand side of (4.5) we get
\begin{align}\label{}
    J(u_{n},v_{n})-J((w,z)_{\delta}) =& (1-\xi_{n}^{+}(\delta ))\langle J'((w,z)_{\delta}),(u_{n},v_{n})\rangle\notag\\
                                     &+\delta\xi_{n}^{+}(\delta )\langle J'((w,z)_{\delta}),
     \frac {J'(u_{n},v_{n})}{\|J'(u_{n},v_{n})\|}\rangle\notag\\
                                     &+o(\|(w,z)_{\delta}-(u_{n},v_{n})\|).
\end{align}
\par Combining~(4.6)~with (4.5) we obtain
     $$~~~~(1-\xi_{n}^{+}(\delta ))\langle J'((w,z)_{\delta}),(u_{n},v_{n})\rangle+\delta\xi_{n}^{+}(\delta )\langle J'((w,z)_{\delta}),
     \dfrac {J'(u_{n},v_{n})}{\|J'(u_{n},v_{n})\|}\rangle$$
     $$\leq o(\|(w,z)_{\delta}-(u_{n},v_{n})\|)+\dfrac{1}{n}\|(w,z)_{\delta}-(u_{n},v_{n})\|.~~~~~~~~~~~~~~~~~~~~~~~~$$
\par
Divide~the above inequality by $\delta$ for~$\delta\neq0$~and
let~$\delta\rightarrow 0$, then we get
    $$ -(\xi_{n}^{+})'(0)\langle J'(u_{n},v_{n}),(u_{n},v_{n})\rangle+
     \|J'(u_{n},v_{n})\| $$
    $$ \leq (o(1)+\dfrac{1}{n})(1+|(\xi_{n}^{+})'(0)|\cdot\|(u_{n},v_{n})\|).~~~~~~~~~~~$$
\par That is,
\begin{center}
     $\|J'(u_{n},v_{n})\| \leq (o(1)+\dfrac{1}{n})\cdot(1+|(\xi_{n}^{+})'(0)|\cdot\|(u_{n},v_{n})\|).$
\end{center}
\par
By (4.4), we only need to show $|(\xi_{n}^{+})'(0)|$~is
uniformly bounded with respect to $n$. Noticing
that
\begin{center}
    $(\xi_{n}^{+})'(0)=\langle(\xi^{+})'(0,0),\dfrac {J'(u_{n},v_{n})}{\|J'(u_{n},v_{n})\|}\rangle,$
\end{center}
by~(4.4) and the third assertion of Lemma 4.4,  we can get that there exists $C>0$ such that
\begin{center}
    $|(\xi_{n}^{+})'(0)|\leq\dfrac{C}{|\|(u_{n},v_{n})\|^{2}-3A(u_{n},v_{n})|}.$
\end{center}
\par Thus we only need to prove that $|\|(u_{n},v_{n})\|^{2}-3A(u_{n},v_{n})|$~has a positive lower bound.~Assume the contrary, then up to a subsequence,
\begin{equation}\label{}
    \|(u_{n},v_{n})\|^{2}-3A(u_{n},v_{n})=o(1).
\end{equation}
\par Since~$\{(u_{n},v_{n})\}\subset \mathcal{N}^{+}$,
\begin{equation}\label{}
    \|(u_{n},v_{n})\|^{2}-A(u_{n},v_{n})=B(u_{n},v_{n}).
\end{equation}
\par From (4.7) and (4.8) we have,
\begin{equation}\label{}
   B(u_{n},v_{n})=\frac{2}{3} \|(u_{n},v_{n})\|^{2}+o(1).
\end{equation}
\par For fixed~$f(x),g(x)$, since
$\max\{\|f\|_{\frac{4}{3}},\|g\|_{\frac{4}{3}}\}<\Lambda$,~there must exist a small positive~$\tau$~such that
\begin{equation}\label{}
    \max\{\|f\|_{\frac{4}{3}},\|g\|_{\frac{4}{3}}\}<(1-\tau)\Lambda.
\end{equation}
\par By the proof of Lemma~2.1, one knows from the derivation of (2.3) that more accurate inequality will occur once (4.10) holds. That is, one has for~$\|(u,v)\|=1$,
\begin{center}
     $B(u,v)<\dfrac{2}{3}(1-\tau)\sqrt{\dfrac{1}{3A(u,v)}}.$
\end{center}
\par Thus by homogeneity,
\begin{equation}\label{}
    B(u_{n},v_{n})<\frac{2}{3}(1-\tau)\sqrt{\frac{\|(u_{n},v_{n})\|^{2}}{3A(u_{n},v_{n})}}\|(u_{n},v_{n})\|^{2}.
\end{equation}
\par Dividing~(4.11)~by~$\|(u_{n},v_{n})\|^{2}$ and letting $n \rightarrow \infty$, we reach a contradiction from (4.4),~(4.7)~and~(4.9), that
\begin{center}
    ${\dfrac{2}{3}}\leq\dfrac{2}{3}(1-\tau).$
\end{center}
This completes the proof.\hfill{$\square$}\\

On the other hand, we have\\
\\
\textbf{Lemma~4.7}~~\textit{Assume~$f(x),g(x)\in L^{\frac{4}{3}}(\Omega)$, both nonzero, and $\max\{\|f\|_{\frac{4}{3}},\|g\|_{\frac{4}{3}}\}<\Lambda$. Then there exists a sequence~$\{(u_{n},v_{n})\}\subset \mathcal{N}^{-}$,~such that when~$n \rightarrow \infty$, it holds}:
\begin{enumerate}
  \item $J(u_{n},v_{n})\rightarrow \theta^{-};$
  \item $J'(u_{n},v_{n})\rightarrow 0.$
\end{enumerate}
\textbf{Proof.}~~By Lemma~4.3, $\mathcal{N}^{-}$~is closed in~$H$.~We use Ekeland's variational principle on~$\mathcal{N}^{-}$~to obtain a minimizing sequence~$\{(u_{n},v_{n})\}\subset \mathcal{N}^{-}$, such that
\begin{enumerate}
  \item $J(u_{n},v_{n})<\inf\limits_{(u,v)\in \mathcal{N}^{-}}J(u,v)+\dfrac{1}{n};$
  \item $J(w,z)\geq J(u_{n},v_{n}) -\dfrac{1}{n}\|(w-u_{n},z-v_{n})\|$~holds for all~$(w,z)\in \mathcal{N}^{-}$.
\end{enumerate}
\par From Lemma 4.1 and the remark before Lemma 4.3 we obtain estimates similar to (4.4):
\begin{center}\label{}
    $0<\tilde{m}\leq \|(u_{n},v_{n})\| \leq \tilde{M}$,
\end{center}
where~$\tilde{m}$~and~$\tilde{M}$~are positive constants. By Lemma 4.5, the rest of the proof is similar to that of Lemma 4.6,  we omit it.\hfill{$\square$}\\

Since system (1.1) is  subcritical, it is not difficult to obtain compactness conditions for the functional $J$. That is,\\
\\
\textbf{Lemma~4.8}~~\textit{Assume~$f(x),g(x)\in L^{\frac{4}{3}}(\Omega)$, then $J$ satisfies $(PS)$ condition, i.e., for any $c\in \mathbb{R}$,
any sequence $\{u_{n}\}\subset H$ for which $J(u_{n})\rightarrow c,J'(u_{n})\rightarrow 0$ as $n\rightarrow+\infty$ possesses a convergent subsequence}.\\

 We are now able to give\\

  \par
  \textbf{Proof of Theorem 1.1.}~ Firstly let us consider the minimization problem~(3.2).~By Lemma 4.6,~there exists~$\{(u_{n},v_{n})\}\subset \mathcal{N}^{+}$~such that as~$n\rightarrow\infty$
\begin{center}
  $J(u_{n},v_{n})\rightarrow \theta^{+}$,~~~~
   $J'(u_{n},v_{n})\rightarrow 0$.
\end{center}
\par Since $J$ satisfies $(PS)$condition by Lemma 4.8, we find a $(w_{1},z_{1})\in cl~(\mathcal{N}^{+})\subset\mathcal{N}^{+}\cup \{(0,0)\}$ such that $J(w_{1},z_{1})= \theta^{+}$,~$J'(w_{1},z_{1})=0 $. By Lemma 4.2, $J(w_{1},z_{1})<0$. Thus
$(w_{1},z_{1})\neq(0,0)$,~which implies~$(w_{1},z_{1})\in \mathcal{N}^{+}$.~We see~$(w_{1},z_{1})$~is a nontrivial weak solution of~(1.1) by Lemma 3.2.

Furthermore,~if~$f$ and $g$~are both positive,~we show the minimizer can be chosen to be a multiple of $(|w_{1}|,| z_{1}|)$. Indeed,~$\|(w_{1},z_{1})\|=\|(|w_{1}|,| z_{1}|)\|$.~Let
  $$(|w_{0}|, |z_{0}|):=\dfrac{(|w_{1}|,| z_{1}|)}{\|(|w_{1}|,| z_{1}|)\|},\;\;\;(w_{0}, z_{0}):=\dfrac{(w_{1}, z_{1})}{\|(w_{1}, z_{1})\|}.$$
 \par Since $(w_{1},z_{1})\in \mathcal{N}^{+}$,~from the proof of Lemma 2.1 we know
$$B(\dfrac{w_{1}}{\|(w_{1},z_{1})\|},\dfrac{z_{1}}{\|(w_{1},z_{1})\|})>0,$$
 thus $B(|w_{1}|,| z_{1}|)\geq B(w_{1},z_{1})>0$, which yields $B(|w_{0}|, |z_{0}|)>0$. By the proof of Lemma 2.1, there exists $t_{1}>0$ such that
$t_{1}\cdot(|w_{1}|,| z_{1}|)\in \mathcal{N}^{+}$.~Since
 $t_{1}\|(w_{1},z_{1})\|\cdot(|w_{0}|, |z_{0}|)=t_{1}\|(|w_{1}|,|z_{1}|)\|\cdot(|w_{0}|, |z_{0}|)=t_{1}\cdot(|w_{1}|,| z_{1}|)\in \mathcal{N}^{+}$,
 we know~$t_{1}\|(w_{1},z_{1})\|$~is the first stationary point of the fibering map in the direction $(|w_{0}|, |z_{0}|)$.
  \par
  Moreover, $(w_{1},z_{1})\in \mathcal{N}^{+}$ is equivalent to $\|(w_{1},z_{1})\|(w_{0}, z_{0})\in \mathcal{N}^{+}$, so $\|(w_{1},z_{1})\|$ is the first stationary point of the fibering map in the direction $(w_{0}, z_{0})$. Since $B(|w_{0}|,| z_{0}|)\geq B(w_{0},z_{0})>0$ and $A(|w_{0}|,| z_{0}|)=A(w_{0},z_{0})$, we can compare the above two roots of the associated fibering map to infer $t_{1}\|(w_{1},z_{1})\|\geq\|(w_{1},z_{1})\|$. That is
\begin{equation}\label{}
    t_{1}\geq 1.
\end{equation}
Taking account of the graph of the fibering map in direction~$(|w_{1}|,| z_{1}|)$, one has from (4.12) and~the fact $t_{1}\cdot(|w_{1}|,| z_{1}|)\in \mathcal{N}^{+}$ that,
\begin{center}
    $J(t_{1}|w_{1}|,t_{1}| z_{1}|)\leq J(|w_{1}|,| z_{1}|)$.
\end{center}
Thus
\begin{center}
    $\theta ^{+}\leq J(t_{1}|w_{1}|,t_{1}| z_{1}|)\leq J(|w_{1}|,| z_{1}|)\leq J(w_{1}, z_{1})=\theta ^{+}$,
\end{center}
from which we know~$(t_{1}|w_{1}|,t_{1}| z_{1}|)$~solves problem~(3.2).~By Lemma 3.2 we see~$(t_{1}|w_{1}|,t_{1}| z_{1}|)$~is a weak solution of system (1.1).

Now we consider the minimization problem~(3.3). By Lemma 4.7, there exists~$\{(u_{n},v_{n})\}\subset \mathcal{N}^{-}$, such that  as~$n\rightarrow\infty$
\begin{center}
  $J(u_{n},v_{n})\rightarrow \theta^{-}$,~~~~
   $J'(u_{n},v_{n})\rightarrow 0$.
\end{center}
Since~$J$~satisfies~$(PS)$~condition by Lemma 4.8, we find a~$\{(w_{2},z_{2})\}\in cl~(\mathcal{N}^{-})=\mathcal{N}^{-}$
such that $J(w_{2},z_{2})= \theta^{-}$,~$J'(w_{2},z_{2})=0 $, and $(w_{2},z_{2})$ is a nontrivial weak solution of (1.1) by Lemma 3.2.

Furthermore,~if~$f$ and $g$~are both positive,~we show the minimizer for (3.3) can be chosen to be a multiple of $(|w_{2}|,| z_{2}|)$. By the proof of Lemma 2.1, there exists~$t_{2}>0$~such that
$t_{2}\cdot(|w_{2}|,| z_{2}|)\in \mathcal{N}^{-}$. Moreover, one calculates that the two fibering maps in direction $(w_{2}, z_{2})$ and $(|w_{2}|,| z_{2}|)$ has the same turning point denoted $t_{0}=\sqrt{\frac{1}{3A(w_2,z_2)}}=\sqrt{\frac{1}{3A(|w_2|,|z_2|)}}$. Thus~$t_{2}>t_{0}$, and by investigating the graph of the fibering map in direction $(w_{2},z_{2})$~one gets
\begin{center}
    $J(t_{2}w_{2},t_{2} z_{2})\leq J(w_{2}, z_{2})$.
\end{center}
\par Now we have
\begin{center}
    $\theta ^{-}\leq J(t_{2}|w_{2}|,t_{2}| z_{2}|)\leq J(t_{2}w_{2},t_{2} z_{2})\leq J(w_{2}, z_{2})=\theta ^{-}$,
\end{center}
so~$(t_{2}|w_{2}|,t_{2}| z_{2}|)$~solves~(3.3).~By Lemma 3.2 we see~$(t_{2}|w_{2}|,t_{2}| z_{2}|)$~is a weak solution of system (1.1).

We finish the proof by showing $\theta_{+}<\theta_{-}$. In fact, if $(w_{2},z_{2})$   is the minimizer of (3.3) satisfying $B(w_{2},z_{2})\leq0$, then
the associated fibering map has only one stationary point, which implies $\theta_{-}>0$ by the proof of Lemma 2.1.
 So $\theta_{+}<\theta_{-}$ by Lemma 4.2. On the other hand, if $(w_{2},z_{2})$ satisfies $B(w_{2},z_{2})>0$, then
 the associated fibering map $\phi(t)$ has two stationary points: $t_{1},t_{2}(=1)$. Thus from the graph of this
  fibering map, we get immediately that $\theta_{+}\leq \phi(t_{1}) <\phi(t_{2})=\theta_{-}$. \hfill{$\square$}
 \vskip 0.4cm
 {\bf Remark} ~~~ For the definition $A(u,v)$ (2.1), we have by H\"{o}lder inequality $$\int_{\Omega}u^2v^2\leq(\int_{\Omega}u^4)^{\frac{1}{2}}(\int_{\Omega}v^4)^{\frac{1}{2}}=\|u\|_4^2\cdot\|v\|_4^2,$$
 then
 $$A(u,v)\geq\mu_1\|u\|_4^4+\mu_2\|v\|_4^4-2|\beta|\|u\|_4^2\|v\|_4^2\geq 2(\sqrt{\mu_1\mu_2}-|\beta|)\|u\|_4^2\cdot\|v\|_4^2,$$
 and the equality  is satisfied for the second inequality if and only if $\mu_1\|u\|_4^4=\mu_2\|v\|_4^2$.
 Thus as $\|(u,v)\|=1, \beta>-\sqrt{\mu_1\mu_2}$, we have $$A(u,v)>0.$$
 Therefore, as $\beta>-\sqrt{\mu_1\mu_2}$, all the proofs are valid, Theorem 1.1  above is still true.\hfill{$\square$}
\par
{\bf Acknowledgements:} The authors thank the referees for their
careful reading and helpful suggestion.

\end{document}